\documentclass[12pt]{article}         %[erg_circle.tex 24.01.14]
\usepackage{amsfonts,amssymb} %mlb,

\usepackage[margin=25mm, a4paper]{geometry}  % for the journal

\def\beq#1#2{\begin{equation} \label{#1} #2 \end{equation}}
\def\bea#1{\begin{eqnarray*} #1 \end{eqnarray*}} \def\a{\!\!\!&\!\!\!\!&}
\def\beaq#1#2{\label{#1} \begin{eqnarray} #2 \end{eqnarray}}
\def\n{\noindent}      
    
\def\IR{{\mathbb{R}}}  \def\IZ{{\mathbb{Z}}}  
\def\toas#1{\stackrel{#1}{\longrightarrow}}
\def\ep{\varepsilon}  \def\phi{\varphi}   
  \def\mod1{\,({\rm mod\ } 1)\,}
\def\t{\tilde} \def\rho{\varrho} 
\def\Var{\rm Var}

\def\function#1{\left\{\!\!\!\begin{array}{ll} #1 \end{array} \right.}
\def\proof{\smallskip \noindent {\bf Proof. \ }}       %start of proof
\def\blanksquare{\,\,\,$\sqcup\!\!\!\!\sqcap$}         %blank  square
\def\qed{\hfill\blanksquare\linebreak\smallskip\par}   %end of proof

\def\thname{Theorem}     \def\lmname{Lemma}      \def\prname{Proposition}
\def\dfname{Definition}  \def\crname{Corollary}  \def\rmname{Remark}

\newtheorem{theorem}{\thname}%[section]   %Numbering: Theorem--Other section
\newtheorem{lemma}{\lmname}%[section]     %{lemma}[theorem]{Lemma}   section
\newtheorem{proposition}[lemma]{\prname} %lemma
\newtheorem{corollary}[lemma]{\crname}   %lemma

\newtheorem{dftn}{\dfname}[section]
\newtheorem{rmrk}[lemma]{\rmname}
\newenvironment{remark}{\begin{rmrk}\rm}{\end{rmrk}}     %lemma

             \catcode`@=11 \@addtoreset{equation}{section} \catcode`@=12

\makeatletter\def\fps@figure{htbp}\makeatother %figure pos: tbp - standard
\newcommand\mlbscale{1pt} %to change: \renewcommand\mlbscale{1.3pt}
\newif\iffigs\figstrue %\newif\iffigs\figsfalse -- to fake figures
\def\bline(#1,#2)(#3,#4)(#5){\put(#1,#2){\line(#3,#4){#5}}}  %straight line

\def\bfig(#1,#2)#3#4{\begin{figure} \begin{center}
    \framebox{\setlength{\unitlength}{\mlbscale}
       \iffigs \begin{picture}(#1,#2) #3 \end{picture}
       \else \begin{picture}(60,10)(0,0)
                   \put(0,0){\framebox(60,10){Figure}} \end{picture} \fi}
    \end{center} \caption{#4} \end{figure}}
\def\Bfig(#1,#2)#3#4{\begin{figure} \begin{center}
    \setlength{\unitlength}{\mlbscale}
       \iffigs \begin{picture}(#1,#2) #3 \end{picture}
       \else \begin{picture}(60,10)(0,0)
                   \put(0,0){\framebox(60,10){Figure}} \end{picture} \fi
    \end{center} \caption{#4} \end{figure}}
\def\bpic(#1,#2)#3{\setlength{\unitlength}{\mlbscale}
    \begin{picture}(#1,#2) #3 \end{picture}}
 \def\rho{\varrho} 

\def\*#1{#1^*}    \def\0#1{\breve#1}  \def\2#1{\acute#1}
\def\G{\Delta}

\def\?#1{} % comments 

\begin{document}%%%-------------------------------------------------%%%%
\title{Ergodicity of a collective random walk on a circle}
\author{Michael Blank\thanks{
        Russian Academy of Sci., Inst. for
        Information Transmission Problems, ~
        e-mail: blank@iitp.ru}
        \thanks{This research has been partially supported
                by Russian Foundation for Basic Research and ONIT RAS program.}
       }
\date{January 24, 2014} %\today % 
\maketitle

\begin{abstract}%
We discuss conditions for unique ergodicity of a collective random walk on 
a continuous circle. Individual particles in this collective motion perform 
independent (and different in general) random walks conditioned by the 
assumption that the particles cannot overrun each other. Additionally to 
sufficient conditions for the unique ergodicity we discover a new and 
unexpected way for its violation due to excessively large local jumps. 
Necessary and sufficient conditions for the unique ergodicity of the 
deterministic version of this system are obtained as well. 
Technically our approach is based on the interlacing property of 
the spin function which describes states of pairs of particles in 
coupled processes under study. 
\end{abstract}%

%%%%%%%%%%%%%%%%%%%%%%%%%%%%%%%%
\section{Introduction}\label{s:intro}

We consider a  collective random walk of a configuration consisting of 
$n$ particles on a unit continuous circle. Each particle without interactions 
with others performs an independent random walk and the interaction 
between particles consists in the prohibition for particles to overrun each 
other. The $i$-th particle in the configuration at time $t\in\IZ_+$ is 
characterized by the position of its center $x_i^t\in S:=[0,1)$, the 
radius $r_i\ge0$ of the ball (representing the particle), and the 
distribution of jumps $P_i$ (i.e. the particle makes a jump equal 
to a random value $\xi$ distributed according to $P_i$). In general 
our theory covers both positive and negative jumps, but to simplify 
presentation we discuss in the Introduction only the case of nonnegative jumps, 
i.e. $P_i([0,1])=1$, leaving the general case to Section~\ref{s:vel-2signs}.

The collective random walk under consideration is a close relative to 
exclusion processes introduced by Frank Spitzer \cite{Sp} and studied 
in a number of publications. One of the most prominent and detailed 
review of statistical properties of such processes considered on a lattice 
and in continuous time can be found in \cite{Lig} (see further references 
therein and \cite{An,AAV,EFM,Pe} for more recent results). 

We say that a  particle configuration $x^t:=\{x_1^t,\dots,x_n^t\}$ is
{\em admissible} if the open balls corresponding to the particles in 
the configuration do not intersect (see Fig.~\ref{f:tasep-c}),  
i.e. it satisfies the inequality:
$$ x_i^t+r_i +r_{i+1} \le x_{i+1}^t \qquad \forall i. $$
The set of all admissible configurations we denote by $X$. 
Here and in the sequel (if the exception is not explicitly mentioned) 
arithmetic operations with metric elements ($x_i,v_i,r_i$, etc.) 
are taken modulo $1$, the comparison between them is 
performed according to the {\em clockwise order} as elements 
of the unit circle, and the indices are taken modulo $n$, i.e. 
$x_{n+1}^t\equiv x_1^t, x_0^t\equiv x_n^t$. 

%%%%%%%%%%%%%%%%%%%%%%%%%%%%%%%%%%%%%%
\def\particle#1{
     \put(0,0){\circle{30}} \bline(0,0)(0,-1)(25) \put(-3,-32){$x_{#1}^t$}
     \bline(0,20)(1,0)(15) \put(5,22){$r_{#1}$}
     \put(0,35){\vector(1,0){70}} \put(20,39){$v_{#1}^t$}}
\def\sparticle#1{
     \put(0,0){\circle{25}} \bline(0,0)(0,-1)(25) \put(-3,-32){$x_{#1}^t$}
     \bline(0,20)(1,0)(12) \put(5,22){$r_{#1}$}
     \put(0,35){\vector(1,0){40}} \put(15,39){$v_{#1}^t$}}
\Bfig(150,50)
      {\footnotesize{
       \thicklines
       \bline(0,10)(1,0)(150)
       \put(30,10){\particle{i}} \put(105,10){\sparticle{i+1}}
       \thinlines \bline(46,28)(1,0)(48) \put(65,32){$\Delta_i^t$}
       \bezier{50}(30,10)(55,-20)(78,10) \put(54,-5){\vector(1,0){4}} %\put(53,-1){$p$}
       \put(78,10){\circle{30}} \bline(78,10)(0,-1)(25) \put(75,-22){$x_i^{t+1}$}
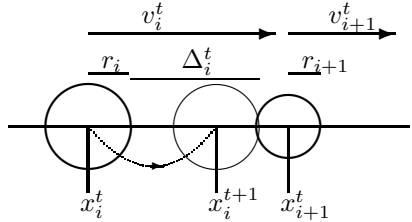
      } } {Exclusion process in continuum \label{f:tasep-c}}
%%%%%%%%%%%%%%%%%%%%%%%%%%%%%%%%%%%%%%

Finally the local dynamics (see Fig.~\ref{f:tasep-c}) of an individual particle 
is defined by the relation
\beq{e:dyn}{x_i^{t+1}:=\min\{x_i^t+u_i^tv_i^t, ~x_{i+1}^t-r_i-r_{i+1}\} ,}%
where the random variable (velocity) $v_i^t$ is chosen according to 
the distribution $P_i$, and the collection of multiplicators $u^t:=\{u_i^t\}_i$ 
with $u_i^t\in\{0,1\}$ represents the updating rule (see below). The random 
variables $\{v_i^t\}_{i,t}$ are assumed to be mutually independent. 

The moment of time when the $i$-th particle 
is stopped by the $i+1$-th particle (i.e. $x_i^t+u_i^tv_i^t>x_{i+1}^t-r_i-r_{i+1}$) 
will be referred as the moment of {\em interaction} between these particles.  
Note that a homogeneous version of systems of this type (when $r_i$ and 
$P_i$ do not depend on $i$) was introduced and studied in \cite{Bl10}. 

Depending on the updating rule $u^t$ discrete time processes under consideration 
may be classified into two types: with {\em parallel} and {\em sequential} updating. 
In the former case all particles are trying to move simultaneously which leads 
to an arbitrary number of simultaneous interactions. In the later case at each 
moment of time only one particle is chosen to jump according to a certain 
rule (e.g. by a random choice) and thus at most a single interaction may 
take place. In terms of $u^t$ the parallel updating means that 
$u_i^t\equiv1$ for all $i,t$. The sequential updating may be realized in 
a number of ways and we shall consider the following two scenarios: \\
(a) random sequential updating: at time $t$ the only positive entry in $u^t$ is 
chosen at random according to a given distribution 
$q:=\{q_1,\dots,q_n\}$ with $\prod_iq_i>0$. \\  % ??? q_i=1/n
(b) deterministic sequential updating: we start by choosing a certain index 
$i$ starting from which the particles are updated clockwise one at a time 
until we reach $i$. Then we repeat the procedure, etc.

If the type of the sequential updating is not specified explicitly we mean that 
it is either (a) or (b). Surprisingly conditions leading to the unique ergodicity in 
both these cases coincide. 
The sequential updating in a sense is equivalent to  continuous time 
collective random walk in which a random alarm clock is attached to each particle 
and the particle moves only when the clock rings. 

The local dynamics (\ref{e:dyn})  
together with a specific updating rule define a finite dimensional Markov chain 
in the phase space of admissible configurations. This local dynamics uniquely 
defines the dynamics of {\em gaps}  
$$ \G_i^t:=(x_{i+1}^t - r_{i+1}) - (x_i^t + r_i) $$
between the particles. Naturally $\G^t:=\{\G_i^t\}_i$ is {\em admissible} if 
$\G_i^t\ge0$ and $\sum_i(\G_i^t+2r_i)=1$.

In terms of gaps the local dynamics of particles may be rewritten as 
$$ x_i^{t+1}:=x_i^t+\min\{u_i^tv_i^t, ~\G_i^t\} ,$$
while the actual dynamics of gaps is described by the relation %
\beq{e:gaps}{ \G_i^{t+1}:=\G_i^t - \min\{u_i^tv_i^t, ~\G_i^t\} 
                                                  + \min\{u_{i+1}^tv_{i+1}^t, ~\G_{i+1}^t\} .}
This shows that the dynamics of gaps $\G^t$ is a Markov chain as well. 

Standard arguments about the compactness of the phase space and the 
continuity of the corresponding Markov operators imply the existence of 
invariant measures of the processes under study. Question about the 
uniqueness of the invariant measure is much more delicate. Our main 
results (Theorems~\ref{t:erg},\ref{t:erg-par-2}) give sufficient conditions 
for the uniqueness of the invariant measure (i.e. for the unique ergodicity) 
in the true random setting when the distributions $P_i$ are nontrivial. 
An unexpected counterexample in the case of parallel updating related to 
excessively large local velocities is constructed in Proposition~\ref{p:non-erg}.  
Despite very week assumptions made in these theorems they cannot 
be applied in the pure deterministic setting when each 
distribution $P_i$ is supported by a single constant local velocity $v_i>0$. 
Nevertheless we show that this deterministic process still might be uniquely 
ergodic (albeit due to different reasons) and give necessary and 
sufficient conditions for this (Theorem~\ref{t:det-erg}). Comparing to our 
earlier note \cite{Bl12}, sufficient conditions for the unique ergodicity are 
formulated in very different terms and became much weaker especially in 
the general non totally asymmetric case. Technically the main improvement 
is that instead of using the principle of ``isolated interactions'' the present 
approach is based on the interlacing property of the spin function 
(see Section~\ref{t:synh}, Lemma~\ref{l:var-local}), which describes 
states of pairs of particles in statically coupled processes under study. 

Let us discuss possible obstacles for the unique ergodicity in the simplest 
setting. Assume for a moment that instead of a system on the continuous 
circle we deal with a finite discrete time lattice system with $L$ sites and 
periodic boundary conditions, i.e. $x_i^t,v_i^t\in\frac1L\IZ,~r_i=1/(2L)$. 
Assume also that each particle jumps with probability $p\in(0,1)$ to its 
nearest right neighboring site if it is not occupied or stays put otherwise. 
In other words we consider the simple discrete time totally asymmetric 
exclusion process with parallel updates. This Markov chain has an important 
property that the probability to reach one state from another in finite time 
is positive, which implies unique ergodicity of this process. On the other 
hand, a simple modification of this process allowing longer particle 
jumps, e.g. by 2 sites instead of one with even $L$, breaks down this 
property. Nevertheless as we shall show in this case and in much more 
complicated case of the continuous circle simple and not especially 
restrictive assumptions on the particle jumps guarantee the ergodicity 
of the dynamics of gaps. In some situations (see below) conditions for 
the unique ergodicity of the original system and of the dynamics of gaps 
coincide but typically this is not the case.

The paper is organized as follows. In Section~\ref{s:basic} we formulate 
our main results. Section~\ref{s:steps} is dedicated to technical 
constructions allowing the analysis of synchronization type phenomena 
and proofs of the main results in the totally asymmetric case. An alternative 
construction based on the dynamics of gaps only is discussed here as well.  
In Section~\ref{s:vel-2signs} we deal with a more general situation 
when the local velocities are taking both signs (and hence particles may move 
toward each other), discuss briefly the lattice version of the collective 
random walk and a more strict version of the particles conflicts resolution.  
The last Section is dedicated to the analysis how unique ergodicity may 
take place in the pure deterministic collective walk.

The author is grateful to the referees for a number of helpful suggestions 
for improvement in the article.

\section{Main results}\label{s:basic}

We start with the totally asymmetric setting, i.e. $v_i^t\ge0$ for all $i,t$. 

\begin{theorem}(Ergodicity)\label{t:erg}
Let the distributions of jumps $v_i^t$ satisfy the following non-degeneracy 
condition %
\beq{e:non-deg}{P(v_i^t > v_j^t)>0 ~~\forall i\ne j, ~\forall t.}
Then the Markov process $\G^t:=\{\G_i^t\}$ with the sequential 
updating is uniquely ergodic, namely for each admissible initial 
$\G^0$ the distributions of the random variables $\G^t$ converges 
as $t\to\infty$ in Cesaro means to a limit which does not depend on 
$\G^0$. In the case of the parallel updating the same claim holds 
if additionally %
\beq{e:vel-small}{P(v_{i}^t < \ep)>0  ~~\forall\ep>0, ~\forall i,t.}
\end{theorem}

The non-degeneracy condition (\ref{e:non-deg}) is a rather distant 
generalization of the simplest law $P_i:=(1-p)\delta_0+p\delta_\sigma$, 
when a particle makes a jump of length $\sigma>0$ with the probability 
$p>0$ or stays put otherwise. The condition~(\ref{e:vel-small}) implies 
that in the parallel updating scenario not all particles interact all the time.  
See also the discussion of the necessity of this condition after the proof 
of Theorem~\ref{t:erg} in Section~\ref{s:steps}. 
The following result shows that the presence of non-interacting particles 
is ``almost'' necessary for the unique ergodicity under parallel updating.

\begin{proposition}\label{p:non-erg}
Let there exist a collection of positive values $v^0:=\{v_i^0\}$ such that 
$\Var(v^0):=\sum_i|v_i^0-v_{i+1}^0|>0$ and %
\beq{e:small-det}{\sum_i v_{i}^0 > 1 - 2\sum_i r_i .}
Then the Markov process $\G^t$ with parallel updating and 
$P(v_i^t\ge v_i^0)=1$ for all $i,t$ has infinitely many ergodic invariant 
measures. 
\end{proposition}  

\begin{remark}\label{r:erg} % ???
Using the static coupling construction developed in Section~\ref{s:steps} one 
can show that each two $x^t$-invariant measures coincide up to a spatial shift, 
which together with the results of Theorem~\ref{t:erg} implies the unique 
ergodicity of the original process under factorization upon spatial shifts.
\end{remark}

Now we formulate sufficient conditions of the unique ergodicity in the general 
non totally asymmetric case when the local velocities may take both 
positive and negative values. 

\begin{theorem}\label{t:erg-par-2}
Let the distributions of jumps $v_i^t$ satisfy one of the following non-degeneracy 
conditions 
\beq{e:non-deg-2}{P(v_i^t > v_j^t)>0 \quad  \forall i\ne j,~\forall t,}
\beq{e:non-deg-3}{ P(v_i^t < v_j^t)>0 \quad  \forall i\ne j,~\forall t.}
Then the Markov process $\G^t:=\{\G_i^t\}$ with sequential updating is 
uniquely ergodic. The same claim holds true in the case of parallel updating 
if additionally at least for one index $i$ we have 
\beq{e:vel-small2}{P(v_i^t v_{i+1}^t<0)>0 .}
\end{theorem}

Observe that if the condition~(\ref{e:vel-small2}) is violated we are back to 
the totally asymmetric case.

Very weak sufficient conditions of the unique ergodicity formulated in 
Theorems~\ref{t:erg},\ref{t:erg-par-2} always include a version of a 
non-degeneracy assumption which gives an impression that in the 
deterministic setting (when the jump distributions $P_i$ are concentrated 
at a single point) the unique ergodicity is excluded. The following result 
addresses this question and shows that the corresponding deterministic 
dynamical system may posses a single nontrivial invariant measure. 

\begin{theorem}\label{t:det-erg} Let $P(v_i^t=v_i)=1$ for some constant 
positive local velocities $v:=\{v_i\}_{i=1}^n$ and $u_i^t\equiv1$. 
Denote $v_{\rm min}:=\min_i\{v_i\}$ and 
$\alpha:=1-2\sum_{i=1}^nr_i -nv_{\rm min}$. 
Then the dynamical system defined by the relation (\ref{e:dyn}) is uniquely 
ergodic if and only if \\
(a) $v_{\rm min}$ is achieved at a single index 
$i_{\rm min}$, \\ 
and one of the following assumptions holds true: \\
(b) $\alpha \ge 0$ and $v_{\rm min}$ is irrational, \\
(c) $\alpha < 0$ and $\alpha+v_{\rm min}$ is irrational. 
\end{theorem}

The proofs of Theorems~\ref{t:erg},\ref{t:erg-par-2} use a 
technical result which (especially due to its deterministic nature) 
is of interest by itself. Assume that the local particle velocities 
$v_i^t$ are given for all $i,t$. 
We say that the particle process satisfies the {\em chain-interacting} 
property if for any initial configuration $x^0$ subsequent particles 
will interact in finite time (see a more detailed version of this 
assumption in Section~\ref{s:steps}) . 

\begin{theorem}(Synchronization)\label{t:synh}
Let the chain-interacting property hold. Then \\
(a) in the case of sequential updating for any initial admissible 
configurations $x^0,\2x^0$ the processes $x^t,\2x^t$ are 
getting synchronized with time, namely 
$$\sum_i^n|\G_i^t-\2\G_i^t|\toas{t\to\infty}0.$$
(b) in the case of parallel updating the same claim holds 
if additionally for infinitely many moments of time $t$ 
there exists $j=j(t)$ such that either $v_j^t<0<v_{j+1}^t$ 
or $v_j^t=0$. \\ %
The chain-interacting property holds in the case of sequential 
updating if
\beq{e:erg-seq-det}{\inf_{i,t}\{v_i^t\} > \frac1n(1-2\sum_i r_i) > 0 ,}
while in the case of parallel updating it is enough to assume that %
\beq{e:desyn}{\forall t_0\ge0,~\forall i\in\{1,\dots,n\}~~\exists\tau_i<\infty: 
                       \quad \sum_{t=t_0}^{t_0+\tau_i}(v_i^t - v_{i+1}^t) > 1 .}
\end{theorem}

%\beq{e:v-small}{ \sum_i v_{i}^t < 1 - 2\sum_i r_i .} %

\section{Synchronization phenomenon and proofs in the totally asymmetric setting}
\label{s:steps}

Define a {\em static coupling} between the processes $x^t,\2x^t$ satisfying the 
relation  (\ref{e:dyn}) as a pair-process $(x,\2x)^t$ in which all random choices 
related to particles with equal indices coincide, i.e. $v_i^t=\2v_i^t$ for all $i,t$. 
In the statically coupled processes we say that the $i$-th pair of particles {\em interacts} 
with the $(i+1)$-th one if at least one of the $i$-th particles in these processes 
is doing so.

\subsection{Construction in terms of particle's positions}

Consider a ``lifting'' of the process $x^t$ acting on the circle $S$ to 
the real line $\IR$ defined by the relation %
\beq{e:lift}{ R(x,t,i):=\function{0 &\mbox{if } t=0, i=1 \\ 
                                \sum_{j=1}^{i-1}\G_j^0 &\mbox{if } t=0, i>1 \\ 
                                 R(x,t-1,i) + \min\{u_i^{t-1}v_i^{t-1},\G_i^{t-1}\} 
                                                                     &\mbox{if } t>0.} }%
Rewriting the definition of a gap in terms of the lifting map we have 
$$ \G_i^t = \function{ R(x,t,i+1) - R(x,t,i)  &\mbox{if } 1\le i<n \\ 
                                 R(x,t,1)+1 - R(x,t,n)  &\mbox{if } i=n } ,$$
and the total distance covered by the $i$-th particle during the time 
from 0 to $t$ is equal to $R(x,t,i)-R(x,0,i)$. 

Let the processes $x^t,\2x^t$ be statically coupled. For each $i$ we 
associate to the $i$-th particle a new random variable 
$$ s_i^t := R(x,t,i) - R(\2x,t,i) ,$$ 
to which we refer as a {\em spin}. It is useful to think about the pair 
of points $R(x,t,i), R(\2x,t,i)$ as a {\em dumbbell} whose disks centers 
lye on two parallel straight lines. Then the spin $s_i^t$ describes the 
state of this dumbbell. 

As we shall see, after the interaction between the $i$-th and $i+1$-th 
particles in one of the processes $x^t$ and $\2x^t$ (or in both of them) 
the corresponding spins become closer to each other in comparison to the 
situation just before the interaction (see Fig.~\ref{f:interaction},\ref{f:spin}). 
On the other hand, this might lead to the increase of the distinction with the 
spins of neighboring $(i-1)$-th and $(i+1)$-th particles, i.e. either 
$|s_{i-1}^t-s_{i}^t|$ or $|s_{i}^t-s_{i+1}^t|$ may grow with time $t$. 
Nevertheless in the worst case the amount to which one of the distinctions 
was enlarged cannot be greater than the amount to which another distinction 
became smaller. The idea of our approach is to show the under dynamics 
the {\em variation} of the collection of spins $s^t$, defined as 
$$ \Var(s^t):=\sum_{i=1}^n|s_i^t-s_{i+1}^t| ,$$ 
is a non-increasing function of the variable $t$ and converges to zero 
monotonically with time. 

\begin{lemma} (Interlacing) \label{l:var-local} 
Let at time $t$ one (or both) of the $i$-th particles in the processes 
$x^t,\2x^t$ interact with the $(i+1)$-th one. Then the interlacing property %
\beq{e:var-l}{\min\{s_i^t,s_{i+1}^t\} \le s_i^{t+1} \le \max\{s_i^t,s_{i+1}^t\} } 
takes place. Additionally, if $s_i^t\ne s_{i+1}^t$ then %
\beq{e:var-s}{|s_i^t-s_{i+1}^t| > |s_i^{t+1}-s_{i+1}^t| .}
\end{lemma}
 \proof Consider the interaction at time $t$ of the $i$-th pair of particles with 
the $(i+1)$-th pair and assume that %
\beq{e:part_int}{ \G_i^t < v_i^t \le \2\G_i^t .} %
This situation is depicted in Fig.~\ref{f:interaction}. 
By definition of the spin, (\ref{e:part_int}) implies that $s_i^t > s_{i+1}^t$. 

%%%%%%%%%%%%%%%%%%%%%%%%%%%%%%
%% Interactions
\Bfig(300,80)
      {\footnotesize{\thicklines
       \put(0,0) {
          \put(0,10){\circle{5}} \put(-8,-3){$\2x_{i-1}^t$} \bline(0,10)(0,1)(50)
          \put(0,60){\circle{5}} \put(-8,70){$x_{i-1}^t$}
          \put(30,10){\circle{5}} \put(25,-3){$\2x_{i}^t$} \bline(30,10)(1,2)(25)
          \put(55,60){\circle{5}} \put(50,70){$x_{i}^t$}
          \put(100,10){\circle{5}} \put(92,-3){$\2x_{i+1}^t$} \bline(100,10)(-1,2)(25)
          \put(75,60){\circle{5}} \put(70,70){$x_{i+1}^t$}
                     }
       \put(140,30){$\Longrightarrow$}
       \put(190,0){
          \put(0,10){\circle{5}} \put(-8,-3){$\2x_{i-1}^{t+1}$} \bline(0,10)(0,1)(50)
          \put(0,60){\circle{5}} \put(-8,70){$x_{i-1}^{t+1}$}
          \put(82,10){\circle{5}} \put(70,-3){$\2x_{i}^{t+1}$} \bline(70,60)(1,-4)(12.5)
          \put(70,60){\circle{5}} \put(53,70){$x_{i}^{t+1}$}
          \put(100,10){\circle{5}} \put(92,-3){$\2x_{i+1}^{t+1}$} \bline(100,10)(-1,2)(25)
          \put(75,60){\circle{5}} \put(78,70){$x_{i+1}^{t+1}$}
             } }}
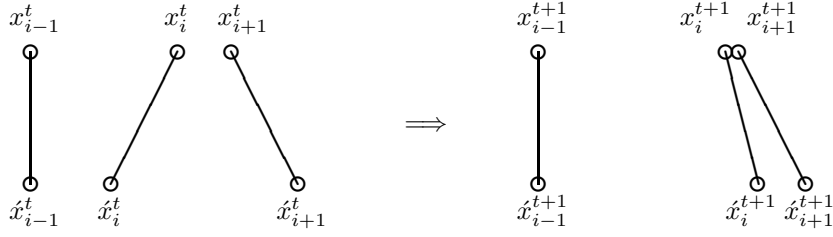
{Isolated interaction of the $i$-th pair with the $(i+1)$-th one. Here we use $x_i^t$ 
  instead of $R(x,t,i)$ to simplify the presentation.
\label{f:interaction}}
%%%%%%%%%%%%%%%%%%%%%%%%%%%%%%

We have $s_i^t=R(x,t,i) - R(\2x,t,i)$ and %
\bea{ R(x,t+1,i) \a= \function{ R(x,t,i+1) - r_i - r_{i+1} 
                                         < R(x,t,i) + v_i^t   &\mbox{if } 1\le i<n \\
                                             R(x,t,1)+1 - r_n - r_{1} 
                                         < R(x,t,n) + v_n^t   &\mbox{if } i=n } \\
      R(\2x,t+1,i) \a= R(\2x,t,i) + v_i^t \le R(\2x,t,i+1) - r_i - r_{i+1} .} %
Therefore %
$$ \min\{s_i^t,s_{i+1}^t\} = s_{i+1}^t \le s_i^{t+1} 
                            < s_i^t = \max\{s_i^t,s_{i+1}^t\} ,$$
which additionally implies (\ref{e:var-s}).

The situation  $\G_i^t \ge v_i^t > \2\G_i^t$ is considered similarly 
exchanging the roles played by the processes $x^t$ and $\2x^t$. 
It remains to study the case 
$$ v_i^t > \max\{\G_i^t,\2\G_i^t \} .$$
This inequality means that after the interaction $s_i^{t+1}=s_i^t$,  
which implies (\ref{e:var-l}), and if additionally $s_i^t\ne s_{i+1}^t$ 
we get (\ref{e:var-s}) as well. \qed

\begin{lemma}\label{l:mon} 
The variation of the spin function $\Var(s^t)$ does not increase under dynamics. 
\end{lemma}

%%%%%%%%%%%%%%%%%%%%%%%%%%%%%%
%% Dynamics of spins
\Bfig(150,130)
      {\footnotesize{\thinlines
        \bline(0,-3)(0,1)(133) %\bline(-10,30)(1,0)(160) 
        \put(-18,122){$s_i^t$}   \put(-25,125){\circle{5}}
        \put(-18,102){$s_i^{t+1}$}  \put(-25,105){\circle*{5}}
        \put(0,3){\circle{5}} %\put(-8,21){$i_k$} 
        \bline(0,3)(-2,1)(10) \bline(0,3)(2,1)(30)
        \bline(30,-3)(0,1)(133)  \put(30,18){\circle{5}}
        \bline(30,18)(1,3)(30)   \put(60,108){\circle{5}}
        \bline(60,108)(2,1)(30) 
        \bline(60,-3)(0,1)(133)  \put(90,122){\circle{5}}
        \bline(90,-4)(0,1)(133)  %\put(79,21){$i_{k+1}$}  
        \bline(90,122)(1,-3)(30) \put(120,32){\circle{5}}
        \bline(120,-3)(0,1)(133) 
        \bline(120,32)(3,-2)(30)   \put(150,12){\circle{5}}
        \bline(150,-3)(0,1)(133) \bline(150,12)(1,1)(10)
        \put(0,10){\circle*{5}}  \bline(0,10)(1,3)(30) \bline(0,10)(-5,1)(10) 
        \put(30,101){\circle*{5}}  \bline(30,101)(2,1)(30) 
        \put(60,115){\circle*{5}}  \bline(60,115)(1,-2)(30) 
        \put(90,55){\circle*{5}}  \bline(90,55)(4,-5)(30) 
        \put(120,18){\circle*{5}}  \bline(120,18)(1,0)(30) 
        \put(150,18){\circle*{5}} \bline(150,18)(2,1)(10)
              }}
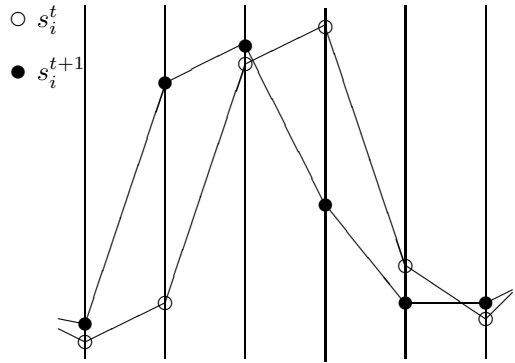
{Change of the spin function under dynamics. $s_i^t$/$s_i^{t+1}$ are 
marked by open/closed circles. \label{f:spin}}
%%%%%%%%%%%%%%%%%%%%%%%%%%%%%%

\proof We say that $[j,k]:=\{j,j+1,\dots,k\}$ is the interval of 
{\em positive monotonicity} of the collection $s:=\{s_i\}$ considered 
as a function of the variable $i$ if $s_j\le s_{j+1}\le\dots\le s_k$, 
and the interval of {\em negative monotonicity} if 
$s_j\ge s_{j+1}\ge\dots\ge s_k$. 

By the property~(\ref{e:var-l}) for each locally maximal interval of 
positive monotonicity $[i,j]$ of the function $s^{t+1}$ we have %
\beq{e:+mon}{\min\{s_i^{t},s_{i+1}^{t}\} \le s_i^{t+1} 
                  \le s_j^{t+1} \le \max\{s_j^{t},s_{j+1}^{t}\} .}
Similarly for each locally maximal interval of negave monotonicity $[i,j]$ %
\beq{e:-mon}{\max\{s_i^{t},s_{i+1}^{t}\} \ge s_i^{t+1} 
                  \ge s_j^{t+1} \ge \min\{s_j^{t},s_{j+1}^{t}\} .}

Consider two consecutive locally maximal intervals of positive monotonicity 
$[i_k^+,j_k^+]$ and $[i_{k+1}^+,j_{k+1}^+]$ of the function $s^{t+1}$. 
Then by~(\ref{e:+mon}) %
\beq{e:bound}{{\rm ind}\max\{s_{j_k^+}^t, s_{j_k^+ +1}^t\}
                  \le {\rm ind}\min\{s_{i_{k+1}^+ +1}^t, s_{i_{k+1}^+}^t\} ,}
where 
$$ {\rm ind}\max\{s_i,s_j\}:=\function{i  &\mbox{if }  s_i\ge s_j \\
                                                         j  &\mbox{otherwise } }, \quad
     {\rm ind}\min\{s_i,s_j\}:=\function{i  &\mbox{if }  s_i\le s_j \\
                                                         j  &\mbox{otherwise } } . $$

The main difficulty in the analysis of the change of variation of the spin 
function is that the intervals of monotonicity of the functions $s^t$ and 
$s^{t+1}$ need not coincide and even may be very different from each 
other (see, e.g. Fig~\ref{f:spin}). Additionally some individual slopes of 
the function $s^{t+1}$ might be much larger than the corresponding 
slopes of the function $s^{t}$. 

Let  $[i_k^+,j_k^+]$ and $[i_k^-,j_k^-]$ be locally maximal intervals of 
positive and negative monotonicity of the function $s$ respectively, and 
let $[i_k,j_k]$ be any collection of non-intersecting intervals. 
(We say that integer intervals do not intersect if they have at most one 
common point.) Then by the triangle inequality %
\beaq{e:var}{ \Var(s)   \nonumber
                \a= \sum_k \left((-s_{i_k^+} + s_{j_k^+}) + (s_{i_k^-} - s_{j_k^-})\right) \\
                \a= 2\sum_k (-s_{i_k^+} + s_{j_k^+}) 
             \ge 2\sum_k (-s_{i_k} + s_{j_k}) .} %

Therefore, combining (\ref{e:+mon}) and (\ref{e:var}) and using 
(\ref{e:bound}) we obtain %
\bea{ \Var(s^{t+1}) \a= 2\sum_k (-s_{i_k^+}^{t+1} + s_{j_k^+}^{t+1}) \\
         \a\le  2\sum_k \left(-\min\{s_{i_k^+ +1}^{t},s_{i_k^+}^{t}\} 
                                    + \max\{s_{j_k^+}^{t}, s_{j_k^++1}^{t}\}\right) \\
         \a\le \Var(s^{t}) ,}
which gives the desired inequality.  \qed

Since the result of Lemma~\ref{l:mon} seems to have an independent 
interest giving a comparison between the variation of two interlacing 
collections of points $s^t$ and $s^{t+1}$, we describe also a sketch 
of an alternative proof of this result based on the induction on the 
number of points $n$. The base of induction -- the case $n=2$ is 
trivial. Indeed, assume for definiteness that $s_1^t<s_2^t$. Then 
$$ s_1^t \le s_1^{t+1} \le s_2^t, \qquad s_1^t \le s_2^{t+1} \le s_2^t .$$
Therefore 
$\Var(s^t)=2|s_1^t-s_2^t|\ge2|s_1^{t+1}-s_2^{t+1}|=\Var(s^{t+1})$. 

It remains to show that the case of general $n>2$ can be reduced to 
the case of a smaller number of points. There are two possibilities: 
there exists at least one interval of monotonicity $J$ of the function 
$s^{t+1}$ of length greater than 1, or all locally maximal 
intervals of monotonicity of the function $s^{t+1}$ are short 
of length 1. In the first case we may remove one of the particles in 
the middle of the interval $J$. This will preserve the interlacing 
conditions (\ref{e:var-l}), but not change the variation 
of $s^{t+1}$ without the removed point, while the variation of $s^t$ 
may only decrease. Thus we get the reduction to the smaller 
number of particles. In the second case we have only short intervals 
of monotonicity of the function $s^{t+1}$ and the  types 
(increasing or decreasing) of corresponding intervals of $s^t$ 
are either opposite, or between two opposite type pairs of 
intervals there is a single pair of intervals with the same types of 
monotonicity. This situation is depicted in Fig.~\ref{f:spin1}.
The analysis of the case of short monotonicity intervals is 
straightforward if one recalls the inequality (\ref{e:bound}).

%%%%%%%%%%%%%%%%%%%%%%%%%%%%%%
%% Dynamics of spins
\Bfig(150,100)
      {\footnotesize{\thinlines
        \bline(0,0)(0,1)(100) \bline(30,0)(0,1)(100) \bline(60,0)(0,1)(100) 
        \bline(90,0)(0,1)(100) \bline(120,0)(0,1)(100) \bline(150,0)(0,1)(100)
        \put( 0, 5){\circle{5}} \bline(0,5)(-1,1)(10) \bline(0,5)(2,1)(30)
        \put(30,20){\circle{5}} \bline(30,20)(2,3)(30)
        \put(60,65){\circle{5}} \bline(60,65)(1,-2)(30)
        \put(90, 5){\circle{5}}  \bline(90, 5)(1,3)(30)
        \put(120,95){\circle{5}}  \bline(120,95)(1,-2)(30)
        \put(150,36){\circle{5}}  \bline(150,36)(1,1)(10)
        \put( 0,15){\circle*{5}} \bline( 0,15)(1,1)(30) \bline( 0,15)(-3,1)(10)
        \put(30,45){\circle*{5}} \bline(30,45)(1,-1)(30)
        \put(60,15){\circle*{5}} \bline(60,15)(1,2)(30)
        \put(90,76){\circle*{5}} \bline(90,76)(1,-1)(30)
        \put(120,47){\circle*{5}} \bline(120,47)(3,1)(30)
        \put(150,56){\circle*{5}} \bline(150,57)(2,-1)(10)
        \put(-18,90){$s_i^t$}   \put(-25,92){\circle{5}}
        \put(-18,70){$s_i^{t+1}$}  \put(-25,72){\circle*{5}}
%        \put(-8,21){$i_k$} \put(0,3){\circle{5}}
              }}
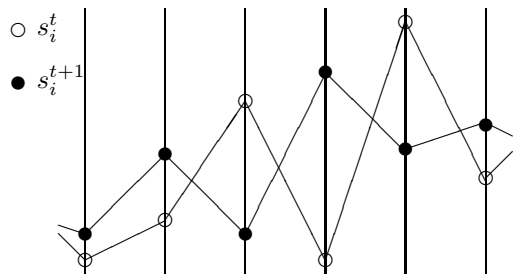
{Oscillation of the spin functions. $s_i^t$/$s_i^{t+1}$ are 
marked by open/closed circles. \label{f:spin1}}
%%%%%%%%%%%%%%%%%%%%%%%%%%%%%%

Now a nontrivial point is to show that under assumptions made above the 
variation $\Var(s^t)$ vanishes with time a.s. To this end one needs to have 
some additional control over particle's interactions. 

We say that the particles with indices $i\le j$ are {\em clockwise chain-interacting} 
if for any initial configuration $x^0$ a.s. there exists a sequence of (random) 
moments of interaction $t_i\le t_{i+1}\le\dots\le t_{j-1}<\infty$ between the 
corresponding particles. In other words, for each $k\in\{0,1,\dots,j-i-1\}$ 
at time $t_{i+k}$  
\beq{e:i1}{v_{i+k}^{t_{i+k}} > \G_{i+k}^{t_{i+k}} . }
Similarly one defines the {\em anti clockwise chain-interaction} when $i>j$. 
If for any pair of indices $i,j$ the (anti) clockwise chain-interacting property holds 
we say that the process satisfies the {\em chain-interacting property}.

\begin{lemma}\label{l:x1}
Let $x^t,\2x^t$ be statically coupled copies of the same particle process with 
sequential updating satisfying the chain-interaction property 
and let $\Var(s^0)>0$. Then 
$\exists\tau=\tau(x^0,\2x^0)<\infty:~~\Var(s^{\tau+1})<\Var(s^\tau)$.
\end{lemma}
\proof By Lemma~\ref{l:var-local} 
$$\Var(s^{t+1})\le\Var(s^t)~~\forall t\ge0$$ 
and we need to show only that for some $t$ this inequality becomes strict. 
Assume from the contrary, that this is not the case, i.e. %
\beq{e:var-eq}{\Var(s^{t+1})=\Var(s^t)~~\forall t\ge0 .} %
By the definition of the sequential updating at the moments of time $t_i$ 
the particle's interactions are isolated in the sense that the particles neighboring 
to the interacting ones make no interactions with other particles. 

Let us show that if there exists an index $i$ such that %
\beq{e:non-mon}{ (s_{i-1}^t-s_i^t)(s_i^t-s_{i+1}^t) < 0 , ~~ 
                            s_i^{t+1}\ne s_i^t , } %
then the variation strictly decreases at time $t$ (i.e. $\Var(s^{t+1})<\Var(s^t)$) . 
The condition (\ref{e:non-mon}) means that $s_i^t$ as a function on $i$ is 
non-monotone at $i$ and changes its value at this index at time $t$. 
This situation is depicted  in Fig.~\ref{f:spin-old}.

%%%%%%%%%%%%%%%%%%%%%%%%%%%%%%
%% Interactions
\Bfig(150,100)
      {\footnotesize{\thicklines
        \bline(0,50)(1,0)(150) \bline(0,0)(0,1)(100)
        \put(20,50){\circle{5}} \put(10,40){$i-1$}
        \bline(20,50)(1,1)(50)
        \put(70,100){\circle{5}} \put(67,40){$i$}
        \bline(70,100)(1,-2)(50)
        \put(120,0){\circle{5}} \put(110,40){$i+1$}
        \put(70,16){\circle*{5}} 
        \thinlines 
        \bline(20,50)(3,-2)(50)  \bline(70,16)(3,-1)(50)
        \put(75,100){$s_i^t$}   \put(75,20){$s_i^{t+1}$}
        \put(10,60){$s_{i-1}^{t}$}  \put(120,10){$s_{i+1}^{t}$}
        \put(70,90){\vector(0,-1){35}}
      }}
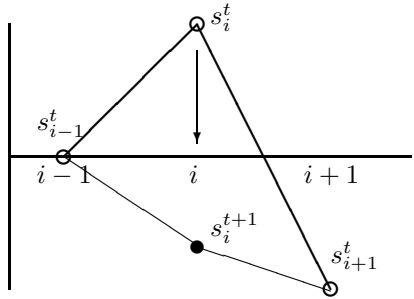
{Change of the $i$-th spin during the isolated interaction. \label{f:spin-old}}
%%%%%%%%%%%%%%%%%%%%%%%%%%%%%%

By the interlacing property~(\ref{e:var-l}) if %
$${\min\{s_{i-1}^t,s_{i+1}^t\} \le s_i^t \le \max\{s_{i-1}^t,s_{i+1}^t\} }$$ %
then $\Var(s^{t+1})=\Var(s^t)$, while the violation of this inequality 
leads by (\ref{e:var-s}) to $\Var(s^{t+1})<\Var(s^t)$. Therefore if 
the index $i$ is the position of a local extremum of the function $s^t$, 
then the variation becomes strictly smaller (see Fig.~\ref{f:spin-old}). 

Now by the chain-interacting property the interaction occurs eventually 
between all neighboring particles and thus the preservation of the variation 
implies that the spin function considered as a function on the index variable 
is monotone (otherwise by the argument above the equality~(\ref{e:var-eq})  
cannot hold). Finally the observation that the spin function is spatially 
periodic (is defined on a circle) shows that this function cannot be monotone, 
unless its variation vanishes (i.e. all spins are equal to each other). 
We came to a contradiction. \qed

\begin{corollary}\label{c:x1}
Let $x^t,\2x^t$ be statically coupled copies of the same particle process with 
sequential updating satisfying  the chain-interaction property. 
Then $\Var(s^t)\toas{t\to\infty}0$. 
\end{corollary}

Indeed, the monotonicity on $t$ (by Lemma~\ref{l:mon}) of the nonnegative 
function $\Var(s^t)$ implies its convergence to a limit, which in turn 
(by Lemma~\ref{l:x1}) cannot differ from zero. 

\begin{lemma}\label{l:x2} $\Var(s^t)\toas{t\to\infty}0$ implies 
$\sum_i^n|\G_i^t-\2\G_i^t|\toas{t\to\infty}0$.
\end{lemma}
\proof By the definition of the spin function 
$$ \G_i^t - \2\G_i^t = s_{i+1}^t - s_i^t ,$$
which implies the claim. \qed

\begin{lemma}\label{l:x3} Let the updating be sequential and let the jump distributions 
satisfy the non-degeneracy condition (\ref{e:non-deg}). Then a.s. in the process
 $x^t$ each pair of particles is either clockwise or anti clockwise chain-interacting. 
\end{lemma}

\proof The condition~(\ref{e:non-deg}) implies that a.s. each particle interacts 
with one (or both) of its neighbors in finite time. Therefore since the total 
number of particles is finite it follows that a.s. after a finite time each 
particle will chain-interact with each other. \qed

\n{\bf Proof of Theorem~\ref{t:erg}}. We start by checking the sequential 
updating case. By Lemma~\ref{l:x3} the chain-interaction property is satisfied. 
Therefore by Corollary~\ref{c:x1} the functional $\Var(s^t)$ vanishes with 
time. Thus the coupling time is finite and hence (see e.g. \cite{BP} for a 
suitable version of the corresponding statement) we get the desired claim.

In the parallel updating case the situation is somewhat more complicated. 
By Lemma~\ref{l:mon} the variation of the spin function cannot increase 
under dynamics and we only need to demonstrate that it strictly decreases 
with positive probability. To this end we make use of the 
condition~(\ref{e:vel-small}), according to which with positive probability 
not all pairs of particles interact simultaneously. 

Assume that at time $t$ the $i$-th pair of particles in the statically 
coupled processes $x^t,\2x^t$ does not interact with the $(i+1)$-th pair. 
Then we can rewrite the parallel updating as $n$ deterministic sequential 
ones with the given local velocities defined in the parallel updating 
starting from the index $i+1$. 

Remark that in the absence of non-interacting particles the reduction to 
the sequential updating cannot be done since otherwise the position of 
the $(i+1)$-th particle will be different at the moment of the sequential 
updating of the $i$-th particle, which will change its position at time $t+1$. 
This is the crucial observation in the proof of Proposition~\ref{p:non-erg} below.

Once we made the reduction to the sequential updating, the results 
proven in that case imply that the variation of the spin function vanishes 
with time and hence we get the unique ergodicity of the corresponding 
gap process for the parallel updating as well. \qed

%{Counter example -- for absence of noninteracting pairs}\bigskip

It is worth noting that in order to guarantee that not all particles 
interact simultaneously all the time it is enough to make an assumption 
\beq{e:vel-small-weak}{P\left(\sum_i v_{i}^t < 1 - 2\sum_i r_i\right)>0 ,}
which is much weaker than (\ref{e:vel-small}). Unfortunately to 
make the simultaneous reduction of both coupled processes from 
parallel to sequential updating we need to find not only a single 
non-interacting particle but a non-interacting dumbbell -- a pair of 
particles which do not interact with their right neighbors. Let us 
show that under the condition~(\ref{e:vel-small-weak}) 
this might not work. Choose a sequence $0<a_0<a_1<\dots<1/2$ 
and let $n=2$, $r_1=r_2=0$. Consider initial configurations of gaps 
$\G_1^0:=a_0, \G_2^0:=1-a_0, \2\G_1^0:=1-a_0, \2\G_2^0:=a_0$. 
To simplify the argument we consider only the deterministic setting, 
choosing the following deterministic sequence of local velocities 
$\{v_1^t,v_2^t\} = \{a_t,a_t\}$. 
Then under dynamics the gaps in the process $x^t$ are equal to 
$\{a_t,1-a_t\}$, while for the process $\2x^t$ they are 
equal to  $\{1-a_t,a_t\}$. Thus the assumption~(\ref{e:vel-small-weak}) 
holds true for any moment of time, only one of the particles in each 
process does not interact with its right neighbor, but in both 1st and 
2nd pairs (dumbbells) one of the particles makes the interaction with 
its right neighbor. 

\bigskip

\n{\bf Proof of  Proposition~\ref{p:non-erg}.} 
Set $a:=\frac1n \left(\sum_i v_i^0 - 1 + 2\sum_i r_i\right)$ and choose 
some $0<b\ll a$. The value $a$ is positive by~(\ref{e:small-det}). 
Consider a configuration of $n$ gaps 
$\t\G:=\{\max\{v_i^0-a/2+b,0 \}\}$. To construct an admissible configuration 
$\G$ we normalize $\t\G$ as follows:
$$ \G:=\left\{ \frac{\t\G_i~ (1 - 2\sum_j r_j)}{\sum_j\t\G_j}  \right\} .$$
By the choice of the parameters $a,b$ for each $i$ we have $v_i^0>\G_i$.

Therefore the application of the dynamics (\ref{e:gaps}) to $\G$ is 
equivalent to the cyclic right shift: $\G_i \to \G_{i+1}$ for all $1\le i<n$ 
and $\G_n\to\G_1$. Thus the configuration $\G$ gives rise to an 
ergodic invariant measure (uniformly distributed on a finite set 
of points) of the process under consideration. 
Noting that choosing different values of the 2nd parameter $b$ we are 
getting different invariant measures we get the claim. \qed

\n{\bf Proof of Theorem~\ref{t:synh} in the case of non negative local velocities}. 
In the sequential updating case the claim about unique ergodicity follows from 
Corollary~\ref{c:x1}, while in the parallel updating case we need additionally 
the condition that for infinitely many moments of time some particles remain 
put to use it instead of the similar probabilistic assumption (\ref{e:vel-small}). 

Therefore we need only to check sufficient conditions for the chain-interacting 
property. In the sequential updating case condition~(\ref{e:erg-seq-det}) 
implies that during at most  $(1-2\sum_i r_i)/\min_i\{v_i^t\}<n$ iterations 
each particle will interact with its nearest right neighbor, which implies the 
property under question. In the parallel updating case the 
condition~(\ref{e:desyn}) plays the same role but does not give explicit 
estimate of the interaction time. \qed

\subsection{Construction in terms of the gap process}\label{s:gaps}

During the discussion of an earlier version of this work a question whether it 
is possible to prove the unique ergodicity using the dynamics of gaps only was posed. 
Here we give a positive answer to this question. Note that despite a certain 
simplification of arguments here we are loosing important information about 
the original particle process and its geometric interpretation. Therefore we prefer 
to discuss both approaches rather than to choose only one of them. 

Recall that the dynamics of gaps is defined by the relation (\ref{e:gaps}). 
Similarly to the particle processes we say that two processes of gaps 
(with the same number of elements)
$\G^t$ and $\2\G^t$ are {\em statically coupled} if 
$v_i^t=\2v_i^t$ for all $i,t$. Define a functional 
$$ V(\G^t,\2\G^t) := \sum_{i=1}^n|\G_i^t-\2\G_i^t| .$$

\begin{lemma}\label{l:gap-dec} 
For a pair of statically coupled processes of gaps $\G^t$ and $\2\G^t$ 
with sequential updating we have %
\beq{e:var-gap}{ V(\G^{t+1},\2\G^{t+1}) \le V(\G^t,\2\G^t) \quad \forall t .}
\end{lemma}

\proof In terms of gaps the interaction between the $i$-th and $(i+1)$-th pair of 
particles at time $t$ in the processes $x^t,\2x^t$ takes place if and only if %
\beq{e:int_gap}{v_i^t>\min\{\G_i^t,\2\G_i^t\}. }
There are 3 possibilities 
\beq{e:part_gap1}{\2\G_i^t \ge v_i^t > \G_i^t,}
\beq{e:part_gap2}{\G_i^t    \ge v_i^t >\2\G_i^t, }
\beq{e:part_gap3}{ v_i^t>\max\{\G_i^t,\2\G_i^t\} .} 
We start with the case (\ref{e:part_gap1}). Then %
\bea{ \G_{i-1}^{t+1}\a=\G_{i-1}^t+\G_i^t, \quad \G_i^{t+1}=0, \\
        \2\G_{i-1}^{t+1}\a=\2\G_{i-1}^t+v_i^t, \quad  \2\G_i^{t+1}=\2\G_i^t-v_i^t.}
Thus %
\bea{ |\G_{i-1}^{t+1}-\2\G_{i-1}^{t+1}| + |\G_{i}^{t+1}-\2\G_{i}^{t+1}| 
  \a= |\G_{i-1}^t+\G_i^t - \2\G_{i-1}^t - v_i^t| + \2\G_i^t-v_i^t \\
  \a=  |(\G_{i-1}^t - \2\G_{i-1}^t) + (\G_i^t - v_i^t)| + \2\G_i^t-v_i^t \\
  \a\le |\G_{i-1}^{t}-\2\G_{i-1}^{t}| + v_i^t - \G_{i}^{t} + \2\G_{i}^{t} - v_i^t \\
  \a= |\G_{i-1}^{t}-\2\G_{i-1}^{t}| - \G_{i}^{t} + \2\G_{i}^{t} \\
 \a\le |\G_{i-1}^{t}-\2\G_{i-1}^{t}| + |\G_{i}^{t}-\2\G_{i}^{t}| .} %
Here the inequality in the 3nd line follows from the triangle inequality. Note 
that this inequality becomes equality if and only if $\G_{i-1}^t \le \2\G_{i-1}^t$ . 
Together with the assumption  (\ref{e:part_gap1}) this implies that two 
consecutive gaps in the process $\G^t$ is less or equal to the corresponding 
gaps in the process $\2\G^t$. 

Similarly in the case  (\ref{e:part_gap2}) we get 
$$ |\G_{i-1}^{t+1}-\2\G_{i-1}^{t+1}| + |\G_{i}^{t+1}-\2\G_{i}^{t+1}| 
   \le |\G_{i-1}^{t}-\2\G_{i-1}^{t}| + |\G_{i}^{t}-\2\G_{i}^{t}| $$
and the inequality takes place if and only if two consecutive gaps in the 
process $\G^t$ is larger or equal to the corresponding 
gaps in the process $\2\G^t$. 

In the remaining case  (\ref{e:part_gap3}) the calculation is even simpler:
$$  \G_{i-1}^{t+1}=\G_{i-1}^t+\G_i^t , \quad \G_i^{t+1}=0,$$
$$ \2\G_{i-1}^{t+1}=\2\G_{i-1}^t+\2\G_i^t, \quad \2\G_i^{t+1}=0. $$
Therefore %
\bea{  |\G_{i-1}^{t+1}-\2\G_{i-1}^{t+1}| + |\G_{i}^{t+1}-\2\G_{i}^{t+1}| 
 \a=  |(\G_{i-1}^t+\G_i^t) - (\2\G_{i-1}^t+\2\G_i^t)| + 0 \\
 \a\le |\G_{i-1}^{t}-\2\G_{i-1}^{t}| + |\G_{i}^{t}-\2\G_{i}^{t}| .}%
Again the equality takes place if and only if two consecutive gaps in the 
process $\G^t$ are both larger or both smaller than the corresponding 
gaps in the process $\2\G^t$.

Observe now that, by the assumption that the process has the sequential updating,  
during the interaction of the $i$-th particle only the $(i-1)$-th and the $i$-th gaps 
may change. Thus the claim (\ref{e:var-gap}) follows. \qed 

Using this result instead of Lemma~\ref{l:var-local} and the functional 
$V(\G^t,\2\G^t)$ instead of the variation of the spin function one can follow 
arguments of the previous Section to prove Theorem~\ref{t:erg}.

\section{Local velocities of both signs and other generalizations}\label{s:vel-2signs} %

\subsection{General (non totally asymmetric) collective random walks}

So far to simplify the setting we assumed that all particles move in the same 
direction, i.e. the local velocities $v_i^t$ have the same (positive) sign. 
The presence of particles moving in opposite directions leads to a
significant modification of the violation of the admissibility condition for 
local velocities. Now one needs to take into account not only the position 
of the succeeding particle, but also its velocity, as well as the corresponding 
quantities related to the preceding particle. In this more general case the 
$i$-th local velocity does not break the admissibility condition if and only if %
\bea{\max\{x^t_{i-1}, x^t_{i-1}+u_{i-1}^tv_{i-1}^t\} + r_{i-1}  %
            \a\le\min\{x^t_i, x^t_i+u_i^tv_i^t\} - r_i \\ %
        \max\{x^t_i, x^t_i+u_i^tv_i^t\} + r_i  
     \a  \le\min\{x^t_{i+1}, x^t_{i+1}+u_{i+1}^tv_{i+1}^t\} - r_{i+1} .}  %
If for some $i\in\{1,2,\dots,n\}$ and $j=i\pm1$ the corresponding
inequality is not satisfied we say that there is a {\em conflict}
between the $i$-th particle and the $j$-th one and one needs to
resolve it. In terms of gaps $\G_i^t$ the inequalities above may be 
rewritten as follows: %
\beq{e:adm-2signs}{\G_{j}^t\ge\max\{u_j^tv_j^t,~-u_{j+1}^tv_{j+1}^t,
                                                      ~u_j^tv_j^t-u_{j+1}^tv_{j+1}^t\},
                   ~~j\in\{i-1,i\} } %

%%%%%%%%%%%%%%%%%%%%%%%%%%%%%%%%%%%%%%
\def\particle#1{
     \put(0,0){\circle{30}} \bline(0,0)(0,-1)(25) \put(-3,-32){$x_{#1}^t$}
     \bline(0,20)(1,0)(15) \put(5,22){$r_{#1}$}
     \put(0,35){\vector(1,0){70}} \put(20,39){$v_{#1}^t$}}
\def\sparticle#1{
     \put(0,0){\circle{25}} \bline(0,0)(0,-1)(25) \put(-3,-32){$x_{#1}^t$}
     \bline(0,20)(1,0)(12) \put(5,22){$r_{#1}$}
     \put(0,30){\vector(-1,0){60}} \put(-25,35){$v_{#1}^t$}}
\Bfig(170,50)
      {\footnotesize{
       \thicklines
       \bline(0,10)(1,0)(170)
       \put(30,10){\particle{i}} \put(150,10){\sparticle{i+1}}    %105
       \thinlines %\bline(46,28)(1,0)(48) %\put(65,32){$\Delta_i^t$}
       \bezier{50}(30,10)(55,-20)(78,10) \put(54,-5){\vector(1,0){4}} %\put(53,-1){$p$}
       \put(78,10){\circle{30}} \bline(78,10)(0,-1)(25) \put(75,-22){$x_i^{t+1}$}
       \put(105,10){\circle{25}} \bline(105,10)(0,-1)(25) \put(103,-22){$x_{i+1}^{t+1}$}
       \bezier{50}(105,10)(127,-20)(150,10) \put(130,-5){\vector(-1,0){4}}
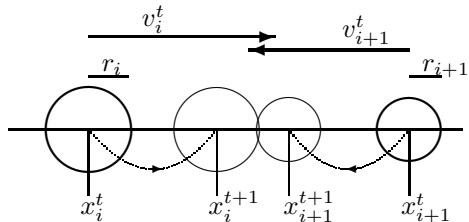
      } } {Velocities of both signs \label{f:vel-sig}}
%%%%%%%%%%%%%%%%%%%%%%%%%%%%%%%%%%%%%%

In distinction to the case of particles moving in the same direction the 
resolution of the conflict between particles is not uniquely defined: to resolve 
a conflict between two mutually conflicting particles moving simultaneously 
in opposite directions (see Fig.~\ref{f:vel-sig}) one needs to specify the 
positions of the particles after the conflict. This can be done in a number of 
ways and we shall consider a {\em natural} resolution of the conflict allowing 
each particle to move with the corresponding velocity as far as possible 
imitating a continuous time motion. Namely in the case of the mutual conflict 
between the $i$-th and the $(i+1)$-th particles, i.e. $v_i^t>0>v_{i+1}^t$, 
the natural resolution of the conflict leads to
$$ x_i^{t+1}:=x_i^t + \frac{\G_i^t v_i^t}{v_i^t-v_{i+1}^t} , \qquad
     x_{i+1}^{t+1}:=x_{i+1}^t + \frac{\G_{i+1}^t v_{i+1}^t}{v_i^t-v_{i+1}^t} . $$

The difference between the condition~(\ref{e:non-deg}) in the formulation 
of Theorem~\ref{t:erg} and the conditions~(\ref{e:non-deg-2}, \ref{e:non-deg-3}) 
in Theorem~\ref{t:erg-par-2} is that the latter are able to deal with particles 
jumping in both directions. 
To show that this is indeed necessary, consider an example with $n=4$ particles 
whose velocity distributions satisfy the relations 
$$P_1((0,1])=P_2([-1,0))=P_3((0,1])=P_4([-1,0))=1 .$$
In this example the particles will eventually meet in pairs (1+2 and 3+4), and 
the pairs will stay put in two random places: 
$$ \G_1^t,\G_3^t\toas{t\to\infty}0, \quad{\rm while}~~\sum_i^4\G_i^t\equiv1 .$$

In the simplest case when the distributions of jumps $P_i$ do not depend 
on the index $i$ it is enough to assume that this common distribution is not 
supported by a single point. 

\bigskip

\n{\bf Proof of Theorem~\ref{t:erg-par-2}}. 
Observe that the constructions developed during the analysis of the 
totally asymmetric setting remain valid in the general case as well. 
The only difference is that additionally to the interaction of a given 
particle with its right nearest neighbor one needs to take into account 
the interactions with the left nearest neighbor when the particle's 
local velocity becomes negative. Fortunately only one of these 
interactions may take place at a given moment of time. 

To apply the machinery developed in Section~\ref{s:steps} we 
need to check that in the case of the mutual conflict the 
interlacing property holds. All definitions made in Section~\ref{s:steps} 
remain valid except for the change in the last line of the 
definition of the lifting~(\ref{e:lift}), where the term 
$\min\{u_i^{t-1}v_i^{t-1},\G_i^{t-1}\}$ should be changed 
to the actual distance covered by the $i$-th particle during its 
jump at time $(t-1)$.

Let the processes $x^t,\2x^t$ be statically 
coupled and the mutual conflict between at least one of the $i$-th 
particles takes place. 
We restrict ourselves only to the situation $\G_i^t\le \2\G_i^t$ (and 
hence $s_i^t\ge s_{i+1}^t$) since the the analysis of the alternative 
situation is completely similar. There are two possibilities.

(a)  $\G_i^t < v_i^t-v_{i+1}^t \le \2\G_i^t$. Then denoting 
$\ell:=\frac{\G_i^t v_i^t}{v_i^t-v_{i+1}^t} < v_i^t$ we get  %
\bea{R(x,t+1,i) \a= R(x,t,i) + \ell, \qquad R(x,t+1,i+1) = R(x,t,i+1) - (\G_i^t-\ell) ,\\
        R(\2x,t+1,i) \a= R(\2x,t,i) + v_i^t, \qquad R(x,t+1,i+1) = R(x,t,i+1) + v_{i+1}^t .}%
Therefore %
\bea{ s_i^{t+1} \a= s_i^t + (\ell-v_i^t) < s_i^t = \max\{s_i^t,s_{i+1}^t\} ,\\
  s_{i+1}^{t+1} \a= s_{i+1}^t + (\G_i^t-\ell+v_{i+1}^t) 
                          > s_{i+1}^t = \min\{s_i^t,s_{i+1}^t\} }%
since $\ell-v_i^t<0$ and  $\G_i^t-\ell+v_{i+1}^t < 0$. 

The observation that $\G_i^{t+1}=0<\2\G_i^{t+1}$ implies 
$s_i^{t+1}>s_{i+1}^{t+1}$, which finishes the analysis of this possibility. 

(b) $v_i^t-v_{i+1}^t > \max\{\G_i^t, \2\G_i^t\}$. Denoting 
$\2\ell:=\frac{\2\G_i^t v_i^t}{v_i^t-v_{i+1}^t} < v_i^t$ we get 
$$ R(x,t+1,i) = R(x,t,i) + \ell, \qquad R(x,t+1,i+1) = R(x,t,i+1) - (\G_i^t-\ell) ,$$
$$ R(\2x,t+1,i) = R(\2x,t,i) + \2\ell, \qquad R(\2x,t+1,i+1) = R(\2x,t,i+1) - (\2\G_i^t-\2\ell) .$$
Therefore 
\bea{   s_i^{t+1} \a= s_i^t + (\ell-\2\ell) \le s_i^t  ,\\
    s_{i+1}^{t+1} \a= s_{i+1}^t - (\G_i^t-\ell) + (\2\G_i^t-\2\ell) \\
                          \a=  s_{i+1}^t - (\G_i^t-\2\G_i^t)(1 - \frac{v_i^t}{v_i^t-v_{i+1}^t})
                          \ge s_{i+1}^t }%
since $\ell\le\2\ell$,  $\G_i^t\le\2\G_i^t$ and  $\frac{v_i^t}{v_i^t-v_{i+1}^t}<1$. 
Eventually we obtain
$$ \min\{s_i^t,s_{i+1}^t\} \le s_{i+1}^{t+1} 
                                         = s_i^{t+1} \le s_i^t = \max\{s_i^t,s_{i+1}^t\} .$$
Additionally a close look to the calculations above shows that if 
$s_i^t\ne s_{i+1}^t$ then %
$$ |s_i^t-s_{i+1}^t| > |s_i^{t+1}-s_{i+1}^{t+1}| ,$$
which gives the analog of the inequality~(\ref{e:var-s}).

It remains to prove that under the assumptions of Theorem~\ref{t:erg-par-2} 
the chain interacting property holds true.
Assume first that the process $x^t$ has sequential updating. Then 
each of the non-degeneracy conditions obviously implies the chain-interacting 
property (in one of the directions). Namely the condition~(\ref{e:non-deg-2}) 
implies the clockwise chain interaction, while the condition~(\ref{e:non-deg-3}) 
implies the anti clockwise chain interaction.
Therefore the claim follows from the same arguments as in the proof 
of the totally asymmetric setting. 

The situation with the parallel updating is slightly more subtle. The point is 
that additionally to the absence of interactions between the nearest particles 
(used in Section~\ref{s:steps}) we get an additional way to make the 
reduction from parallel to consecutive updating. Indeed, if at time $t$ two 
consecutive particles have opposite local velocities then there is at least 
another pair of consecutive particles satisfying this property. Moreover 
among such pairs there is at least one, say $j$ and $j+1$ such that 
$v_j^t<0<v_{j+1}^t$. Therefore these two particles do not interact and  
hence in this situation one can make the reduction to the sequential 
updating, which starts at the index $j+1$ and goes up to the index $j$. 

Using the above trick we make the reduction to the sequential updating 
if the event described in the assumption~(\ref{e:vel-small}) takes place.
The remaining part follows the same arguments as in the proof of 
Theorem~\ref{t:erg}. \qed

\n{\bf Proof of Theorem~\ref{t:synh} in the case of local velocities 
taking both signs}. Additionally to the already proven part related 
to the local velocities of the same sign, we use the condition that 
for infinitely many moments of time $t$ there exists $j=j(t)$ such 
that $v_j^t<0<v_{j+1}^t$ in order to make the reduction from 
the parallel updating to the sequential one at these moments 
of time. After this reduction one applies the same arguments 
as in the probabilistic setting. \qed

\subsection{Strict exclusion}\label{s:strict-exc} %

In all our previous constructions we have considered only those rules resolving particle 
conflicts allowing the particles to move as far as possible according to their local 
velocities. On the other hand, as we already mentioned the collective random walk 
under consideration is a generalization of the simple exclusion process, where a 
particle moves to the neighboring site only if the latter is not occupied by another 
particle. From this point of view it seems to be natural  to consider a similar conflict 
resolution rule. Namely we say that the particle process $x^t$ satisfies the 
{\em strict exclusion} rule if in the case of a conflict the corresponding particle 
stays put, i.e. 
$$ {\rm if} ~~v_i^t>\G_i^t ~~{\rm or}~~ -v_i^t>\G_{i-1}^t ~~{\rm then}~~ 
    x_i^{t+1}=x_i^t .$$
It turns out that this ``natural'' conflict resolution rule typically leads to a non-ergodic 
behavior. 

\begin{proposition}\label{t:strict}  Let $P_i([-\ep,\ep])=0$ for some $\ep>0$ 
and all $i$. Then for $n$ large enough the process $\G^t$ is non-ergodic. 
\end{proposition}

\proof Let $n>1/\ep$. Consider a configuration $x:=\{x_i\}_{i=1}^n$ such that 
$\max_i \G_i<\ep$. Then under the assumptions of the Theorem, due to the strict 
exclusion interaction rule, all particles in the configuration $x$ stay put under 
dynamics. Hence the Dirac measure $\delta_x$ supported by the configuration $x$ 
is invariant under dynamics as well as the Dirac measure supported by the 
sequence of the corresponding gaps $\G_i$. Passing from the configuration $x$ 
to close enough configuration $\2x$ such that $\2\G_i<\ep$ we are getting 
infinitely many different fixed points of the process $x^t$ having different 
configurations of gaps. This proves the non-ergodicity of the process $\G^t$. \qed

In distinction to the non-strict exclusion case here it is much more difficult 
to give sufficient conditions of the unique ergodicity. At present we 
only can formulate the following hypothesis. 

\bigskip
\n{\bf Hypothesis}. Let the assumptions~(\ref{e:non-deg}, \ref{e:vel-small}) 
hold true. Then the process of gaps $\G^t$ with either sequential or parallel 
updating  is uniquely ergodic. 

%\negvskip\negvskip

\subsection{Lattice exclusion process}
Observe that ergodic type results for lattice versions of the problems under consideration 
being nontrivial as well may be derived from the present results. 

In the lattice setting all elements of a configuration $x$ belong to a finite set 
$\{0, 1/n, 2/n,\dots,(n-1)/n\}$ for some $n\in\IZ_+$ which defines the number 
of lattice sites. The radius of a ball representing a particle satisfies the 
condition $nr_i\in\{0,k+1/2\}$ with $k\in\IZ_+$, and the jump distribution 
is supported by the lattice points $P_i(\cup_{j=0}^{n-1}\{j/n\})=1$. 

Despite an apparent significant difference between the behavior of the 
lattice processes in the cases when $r_i\equiv0$ and $r_i>0$ 
(in the former case an arbitrary number of particles may share the 
same lattice site, while in the latter case at most one particle is 
allowed per lattice site) the ergodicity conditions turn out to be 
the same.

\section{Unique ergodicity in the deterministic setting}

Now we address the question of unique ergodicity of the collective 
walk in the deterministic setting. This means that the jump distributions 
are supported by single points: $P(v_i^t=v_i)=1$ for some constant 
local velocities $v:=\{v_i\}_{i=1}^n$. 

We start the analysis from the case when the local velocities $v_i$ take 
both positive and negative values and consider the partition of the 
set if indices into groups of consecutive indices of three types 
corresponding to negative, positive and zero velocities, e.g. the 
configuration of signs ~$++ 000--+--$~ has 5 groups of all 3 types. 
Since there are oppositely signed velocities the number of groups 
is greater than one. 

\begin{theorem}\label{t:det-erg+-} 
The process of gaps $\G^t$ (with either sequential or parallel updating) 
is uniquely ergodic if and only if the number of different groups is 
at most three, and the only group of zero velocity particles (if it exists) 
consists of a single element which is located after the group of positive 
and before negative particles (i.e. $+++~0--$).
\end{theorem}

\proof If the conditions of Theorem are satisfied each initial configuration 
of gaps converges in finite time to the configuration having a single 
positive gap of length $(1-2\sum_i\G_i)$, which proves the unique 
ergodicity. The presence of two zero velocities or more than two 
groups of signed velocities obviously contradicts to the unique 
ergodicity. The observation that the wrong location of the unique 
zero velocity particle (i.e. $++--0$) leads to the presence of the 
invariant measure supported by two points and sensitively depending 
on the position of this particle finishes the proof. \qed

Now we are ready to address a more interesting situation of local 
velocities of the same (say positive) sign. If the updating is sequential,  
Theorem~\ref{t:synh} gives a sufficient condition of the unique ergodicity 
$$ v_{\rm min}:=\min_i\{v_i\} > \frac1n(1-2\sum_i r_i) .$$ 
This condition is probably non optimal and one is tempted to 
weaken it to 
\beq{e:erg-seq-det-w}{\sum_iv_i > 1-2\sum_i r_i .} 
Unfortunately a simple example with two zero velocities 
$v_1=v_2=0,~v_3=v_4=\dots=v_n=1$, which definitely 
satisfies~(\ref{e:erg-seq-det-w}), leads to infinitely many 
invariant measures of the gap process. 

\bigskip

In the case of the parallel updating we able to get both necessary 
and sufficient conditions for the unique ergodicity, formulated in 
Theorem~\ref{t:det-erg}.

\bigskip

\n {\bf Proof of Theorem~\ref{t:det-erg}}. 
Assume that the assumption (a) holds true. 
Since the system is translationally invariant without any loss of generality 
we may assume that the only minimum is achieved at the index 
$i_{\rm min}=n$. The key point to our argument is that for any initial 
admissible particles configuration $x^0$ there exists a finite time 
$t_n=t_n(x^0,v,\{r_i\})$ such that for each $t>t_n$ 
the particle configuration $x^t$ satisfy the property: %
\beq{e:lin}{\G_1^t=\G_2^t=\dots\G_{n-1}^t=\beta, \quad \G_n^t\ge\beta ,}
for a certain $0<\beta\le v_n$. If $\beta<v_n$, then $\G_n^t=\beta$. 

Indeed, if this is the case, starting from the time $t_n$ our dynamical system is 
a direct product of $n$ identical irrational rotation maps 
$$x_i^{t+1}:=x_i^t+\beta~~{\rm mod}(1) .$$
This direct product system possesses a number of invariant measures, but 
the property (\ref{e:lin}) defines a unique invariant measure uniformly 
distributed on the segment given in the coordinates $(x_1,\dots,x_n)$ by the 
relations  
$$ x_1=x_n-(n-1)\beta, ~x_2=x_n-(n-2)\beta, \dots, x_{n-1}=x_n-\beta, $$
provided the number $\beta$ is irrational.

Let us prove that the property (\ref{e:lin}) holds true. Observe that after at most 
$$t_{n-1}:=\left(1-2\sum_{i=1}^nr_i\right)/(v_{n-1}-v_n)$$ 
iterations the $(n-1)$-th particle will catch with the $n$-th one and thus for 
each $t\ge t_{n-1}$ the corresponding gap $\G_{n-1}^t$ is exactly equal to 
the length of jump that the $n$-th particle will perform at time $t$, 
in particular $\G_{n-1}^t\le v_n$. 

Similarly 
$$ \forall t\ge t_{n-2}:=t_{n-1} + \left(1-2\sum_{i=1}^nr_i\right)/(v_{n-2}-v_n)$$ 
the gap $\G_{n-2}^t$ will match the length of jump that the $n$-th particle, etc. 
Eventually after at most
$$ t_{1} := \left(1-2\sum_{i=1}^nr_i\right) 
                \left( \frac1{v_{n-1}-v_n} + \frac1{v_{n-2}-v_n} + \dots 
                                                      + \frac1{v_{1}-v_n}  \right) $$
iterations all particles will start moving synchronously. 

If the assumption (b) holds true, then for for all $t\ge t_1$ 
$$ \G_n^t \ge 1-2\sum_{i=1}^nr_i -(n-1)v_n = \alpha+v_n $$ 
and hence the $n$-th particle will stop interacting with others and will 
perform the pure rotation through the irrational angle $v_n$, which guarantees 
the unique ergodicity. Note that the rationality of the rotation leads to the 
presence of infinitely many invariant measures.

If (c) holds true then $\G_n^t=\alpha+v_n<v_n$ for all $t\ge t_1$ and 
due to the interactions with other particles the $n$-th one will perform the 
pure rotation through another irrational angle $\alpha+v_n$. 

It remains to check that the violation of the assumption (a) implies the 
absence of unique ergodicity. If the minimum is achieved at several 
indices then each of the slowest particles will generate a ``train'' of faster 
particles following it in the same manner as it has been shown for the 
case of the single minimum. If %
\beq{e:small}{1-2\sum_{i=1}^nr_i > nv_{\rm min} }%
then using the same argument as in the case of the single minimum one 
finds a partition of particles into groups following one of the slowest 
particles. In each group the particles are moving at the same speed as 
the leading one. By (\ref{e:small}) each group may be slightly shifted 
not perturbing its motion and the motion of the other particles. Therefore 
the system possesses an infinite number of invariant measures, 
corresponding to the trajectories of the perturbed configurations. 

If the inequality  (\ref{e:small}) is violated, then there are infinitely many 
different particle configurations $x^0$ such that $\G^0<v_{\rm min}$ 
and this property holds for any $t>0$. Thus under dynamics we get 
$\G_i^{t+1}:=\G_{i+1}^t$ for all $i,t$. 
Therefore the obtained configuration is a periodic point of the process $\G^t$, 
which again implies the presence of a infinite number of invariant 
measures. \qed

Note that the ``most natural'' case of identical particles and equal 
constant local velocities, i.e. when 
$ r_1=\dots=r_n, ~v_1=\dots=v_n$, does not satisfy the 
conditions of unique ergodicity.

\end{document}